arXiv:math/0211468v1 [math.NT] 13 Nov 2002

# Cyclotomic Swan subgroups and Primitive Roots

Timothy Kohl

Office of Information Technology, Boston University

Boston, Massachusetts 02215

and

Daniel R Replogle

Department of Mathematics and Computer Science, College of Saint Elizabeth

Morristown, New Jersey 07960

Correspondence address: Daniel R Replogle, Department of Mathematics and Computer Science, College of Saint Elizabeth, 2 Convent Road, Morristown, NJ 07960

Email addresses: tkohl@bu.edu dreplogle@cse.edu

Proposed running head: Swan subgroups and primitive roots

Abstract: Let $K_m = \mathbb{Q}(\zeta_m)$ where $\zeta_m$ is a primitive $m$th root of unity. Let $p > 2$ be prime and let $C_p$ denote the group of order $p$. The ring of algebraic integers of $K_m$ is $\mathcal{O}_m = \mathbb{Z}[\zeta_m]$. Let $\Lambda_{m,p}$ denote the order $\mathcal{O}_m[C_p]$ in the algebra $K_m[C_p]$. Consider the kernel group $D(\Lambda_{m,p})$ and the Swan subgroup $T(\Lambda_{m,p})$. If $(p,m) = 1$ these two subgroups of the class group coincide. Restricting to when there is a rational prime $p$ that is prime in $\mathcal{O}_m$ requires $m = 4$ or $q^n$ where $q > 2$ is prime. For each such $m$, $3 \leq m \leq 100$, we give such a prime, and show that one may compute $T(\Lambda_{m,p})$ as a quotient of the group of units of a finite field. When $h^+_{mp} = 1$ we give exact values for $|T(\Lambda_{m,p})|$, and for other cases we provide an upper bound. We explore the Galois module theoretic implications of these results.

## Section 1: Introduction

For any algebraic number field $N$, we denote its ring of algebraic integers by $\mathcal{O}_N$. Given an abelian, tame (i.e. at most tamely ramified) Galois extension of number fields $L/K$ with Galois group $G$, if $\mathcal{O}_L$ is a free $\mathcal{O}_K[G]$-module, one says $L/K$ has a trivial Galois module structure. Equivalently, one says $L/K$ has a normal integral basis. In [4] the authors show that for any base field $K \neq \mathbb{Q}$ ($\mathbb{Q}$ the field of rational numbers), there will exist fields $L$ so that $L/K$ is of prime degree $p > 2$ so that $L$ does not have a normal integral basis over $K$. See [4, Theorem 2]. The proof, however, does not explicitly exhibit any such primes. This "omission" has been rectified in the cases where $K$ is imaginary quadratic or cyclotomic, [12] and [2]. We continue to pursue the case where $K$ is a cyclotomic field. Throughout, let $C_p$ denote the cyclic group of order $p$, with $p > 2$ a prime.

Let $\Lambda$ denote the order $\mathcal{O}_K[G]$ in the group algebra $K[G]$, where $G$ is a finite group of order $n$. For each tame Galois extension $L/K$ one has $\mathcal{O}_L$ is a locally free rank one $\Lambda$-module, each such extension determines a Galois module class, $[\mathcal{O}_L]$ in the locally free class group $Cl(\Lambda)$. The set of such classes is denoted $R(\Lambda)$. In [9], McCulloh shows, for all abelian groups $G$, that $R(\Lambda)$ is a subgroup of $Cl(\Lambda)$. In [8], McCulloh gives a convenient explicit description of $R(\Lambda)$ in the case $G$ is $p$-elementary abelian. We will make use of a consequence of this description in the case $G \cong C_p$.

In relative Galois module theory, one studies other subgroups of $Cl(\Lambda)$. The easiest one to describe is the kernel group, $D(\Lambda)$. This consists of all classes of locally free $\Lambda$-modules that become trivial upon extension of scalars to the maximal order in $K[G]$. In [16], Ullom studied a subgroup of $D(\Lambda)$ called the Swan subgroup. Let $\Sigma = \sum_{g\in G} g$, and for each $r \in \mathcal{O}_K$ so that $r$ and $n$ are relatively prime, define the Swan module $\langle r, \Sigma\rangle$ by $\langle r, \Sigma\rangle = r\Lambda + \Lambda\Sigma$. Swan modules are locally free rank one $\Lambda$-ideals and hence determine classes in $Cl(\Lambda)$ [16]. The set of classes of Swan modules is $T(\Lambda)$, the Swan subgroup.

Let $\mathcal{O} = \mathcal{O}_K$, $\Gamma = \Lambda/(\Sigma), \overline{\mathcal{O}} = \mathcal{O}/n\mathcal{O}$, $\phi$ and $\overline{\phi}$ denote the canonical quotient maps, $\epsilon$ the augmentation map, and let $\overline{\epsilon}$ be the map induced by $\epsilon$. Consider the fiber product:

$$\begin{array}{ccc} \Lambda & \xrightarrow{\ \phi\ } & \Gamma \\ {\scriptstyle\epsilon}\downarrow & & \downarrow{\scriptstyle\overline{\epsilon}} \\ \mathcal{O} & \xrightarrow{\ \overline{\phi}\ } & \overline{\mathcal{O}} \end{array}$$

When $G$ is abelian (or more generally, when the group algebra $K[G]$ satisfies an "Eichler condition"– see [10] or [3]) this fiber product yields the exact Mayer-Vietoris sequence of Reiner-Ullom:

$$1 \longrightarrow \Lambda^* \longrightarrow \mathcal{O}^* \times \Gamma^* \xrightarrow{h} \overline{\mathcal{O}}^* \xrightarrow{\delta} D(\Lambda) \longrightarrow D(\Gamma) \oplus D(\mathcal{O}) \longrightarrow 0,$$

where for any ring $A$ we denote its group of units by $A^*$. From [16] we have that the image of $s \in \overline{\mathcal{O}}^*$ under $\delta$ is $[s, \Sigma]$, the class of the Swan module $\langle s, \Sigma\rangle$, hence $Im(\delta) = T(\Lambda)$. Therefore, $T(\Lambda)$ is a subgroup of $D(\Lambda)$ and $T(\Lambda) \cong \overline{\mathcal{O}}^*/h(\mathcal{O}^* \times \Gamma^*)$. The map $h$ is given by $(u, v) \mapsto \overline{u} \cdot \overline{v}^{-1} = \overline{\phi}(u)\overline{\epsilon}(v)^{-1}$.

In the case $G$ is $p$-elementary abelian one may establish a relationship between $T(\Lambda)$, $R(\Lambda)$, and $D(\Lambda)$. We state a special case of this relationship here.

**Proposition 1.1 (cf. [4, Proposition 4]).** *Let $G \cong C_p$ and let $T^w(\Lambda)$ denote those classes in $T(\Lambda)$ expressible as $w^{th}$ powers. Then one has $T^{\frac{p-1}{2}}(\Lambda) \subseteq R(\Lambda) \cap D(\Lambda)$.*

From this we have the following immediate corollary, a special case of [4, Corollary 7].

**Corollary 1.2 (cf. [4, Corollary 7]).** *If $G \cong C_p$ and $T(\Lambda)$ is not of exponent dividing $\frac{p-1}{2}$ then $K$ has a tame degree $p$ extension with a nontrivial Galois module structure.*

Henceforth, we assume $G \cong C_p$. A recent area of study (see [2], [6], [11], and [12] for example) has been, for certain base fields $K$, finding explicit $p > 2$ so that one may compute a lower bound on $T(\Lambda)$ and use this to show $R(\Lambda)$ is nontrivial. Finding such an explicit prime $p$ shows there exists a tame Galois extension $L$ of $K$, of degree $p$ so that $\mathcal{O}_L$ is not a free $\Lambda$-module via Corollary 1.2. We express this by saying $L/K$ has a nontrivial Galois module structure for $p$. Specifically, [2] gives an explicit prime $p$ for each cyclotomic field of class number 1 so that there is a degree $p$ extension with nontrivial Galois module structure. For any field not of class number equal to 1, one may show it has a quadratic extension with a nontrivial Galois module structure. Hence, one has a satisfactory result in the case $K$ is a cyclotomic field. Nonetheless, as [2] uses primes that split in $K$, we ask if one may find a rational prime that is prime in $\mathcal{O}_K$ so that $K$ has a tame Galois extension of this order with a nontrivial Galois module structure. In a different direction, [12] uses the fact for $p > 2$ unramified in $K/\mathbb{Q}$ that $T(\Lambda) = D(\Lambda)$ to compute bounds on $R(\Lambda) \cap D(\Lambda)$ when $K$ is an imaginary quadratic field.

We set our notation and outline our approach. Let $K_m = \mathbb{Q}(\zeta_m)$ where $\zeta_m$ is a primitive $m$th root of unity. Let $p > 2$ be prime so that $(p, m) = 1$. The ring of algebraic integers of $K_m$ is $\mathcal{O}_m = \mathbb{Z}[\zeta_m]$. The order $\mathcal{O}_m[C_p]$ in the group algebra $K_m[C_p]$ is denoted $\Lambda_{m,p}$. Section 2 shows for there to be a rational prime $p > 2$ that is prime in $\mathcal{O}_m$ requires $m = 4$, $q^n$, or $2q^n$ for some prime $q > 2$, and that, moreover, one need not consider the case $m = 2q^n$. We find the least prime $p$ with $p > 2$ so that $p$ is prime in $K_m$ when $m = 4$ or $q^n$ with $q > 2$ prime and $3 \leq m < 100$. In Section 3 we give an argument showing $T(\Lambda_{m,p}) = D(\Lambda_{m,p})$ is computable as the cokernel of the map $h$ in the Mayer-Vietoris sequence of Reiner-Ullom. This yields Galois module theoretic consequences and shows that for $p > 2$ prime in $\mathcal{O}_m$ that $T(\Lambda_{m,p})$ is the quotient ring of the units of the finite field $\mathbb{Z}[\zeta_{mp}]/(1 - \zeta_p)\mathbb{Z}[\zeta_{mp}]$. Section 4 is devoted to showing to what extent this is feasibly computable and what tools assist in the computation.

## Section 2: Ramification in Cyclotomic Fields and Inert Primes

A standard result is that a (rational) prime $p$ ramifies in $K_m = \mathbb{Q}(\zeta_m)$ if and only if $p$ divides $m$. So assume $p \nmid m$. Let $f$ be the least such $f$ so that $m$ divides $k = p^f - 1$. From here onward, let $\phi$ denote the Euler $\phi$-function, and set $r = \phi(m)/f$. Then [17, Theorem 2.13] shows $(p)$ factors in $\mathcal{O}_m$ into $r$ prime ideals of residue degree $f$. In particular, if $p \nmid m$ then $p$ is unramified, and if $f = \phi(m)$ then $p$ remains prime in $\mathcal{O}_m$. In the latter case we say $p$ is inert in $K_m/\mathbb{Q}$.

**Proposition 2.1.** *Let $K_m = \mathbb{Q}(\zeta_m)$, then a prime $p > 2$ is inert in $K_m/\mathbb{Q}$ if and*

*only if p is a primitive root mod m. Hence, the only time there is such a prime is if* $m = 4, q^n$*, or* $2q^n$*. Conversely, for each such m, one may find the least such* $p > 2$ *so that* $p \equiv r \bmod m$ *where r is the least primitive root mod m.*

*Proof.* Let $p > 2$ be so that $p \nmid m$. Let $f$ be minimal such that $m$ divides $k = p^f - 1$. and set $r = \phi(m)/f$. Above we noted that $p$ is inert if and only if $r = 1$. In this case $f = \phi(m)$. But if $f = \phi(m)$ then $\phi(m)$ is the minimal value of $f$ so that $p^f \equiv 1 \bmod m$. Therefore, $p$ is a primitive root mod $m$. However, the Primitive Root Theorem of elementary number theory says there is an $r$ so that $r^f \equiv 1 \bmod m$ with $f = \phi(m)$ ($f$ minimal) for $m > 2$ if and only if $m = 4, q^n$, or $2q^n$ for some prime $q > 2$.

In many cases the least primitive root, $r$, is not prime. Observe, however, that if $r$ is the least primitive root mod $m$ then any $s$ so that $r \equiv s \bmod m$ is a primitive root. As any arithmetic progression contains infinitely many primes (Dirichlet), given an $r$ that is a primitive root mod $m$ we may choose a prime $p \equiv r \bmod m$ which is a primitive root mod $m$. Specifically, for each $m$ so that $m = 4$, $q^n$, or $2q^n$ find the least primitive root mod $m$ and denote this by $r_1$. If $r_1 > 2$ and $r_1$ is prime let $p = r_1$. Otherwise consider the arithmetic progression $r_{i+1} = im + r_i$ for $i = 1, 2, 3, \ldots$. Denote the first element of this progression which is prime by $p$. □

Now note for $m \equiv 2 \bmod 4$, one has $K_m = K_{\frac{m}{2}}$. Hence, one normally only considers cyclotomic fields $\mathbb{Q}(\zeta_m)$ with $m$ not congruent to 2 mod 4. Nonetheless, we first present the table one obtains using the algorithm, where $3 \leq m < 100$ with $m = 4$, $q^n$, or $2q^n$.

Table 2.2

| $m$ | $r$ | $p$ | $m$ | $r$ | $p$ |
|---|---|---|---|---|---|
| 3 | 2 | 5 | 4 | 3 | 3 |
| 5 | 2 | 7 | 6 | 5 | 5 |
| 7 | 3 | 3 | 9 | 2 | 11 |
| 10 | 3 | 3 | 11 | 2 | 13 |
| 13 | 2 | 41 | 14 | 3 | 3 |
| 17 | 3 | 3 | 18 | 5 | 5 |
| 19 | 2 | 59 | 22 | 7 | 7 |
| 23 | 5 | 5 | 25 | 2 | 127 |
| 26 | 7 | 7 | 27 | 2 | 29 |
| 29 | 2 | 31 | 31 | 3 | 3 |
| 34 | 3 | 3 | 37 | 2 | 113 |
| 38 | 3 | 3 | 41 | 6 | 47 |
| 43 | 3 | 3 | 46 | 5 | 5 |
| 47 | 5 | 5 | 49 | 3 | 3 |
| 50 | 3 | 3 | 53 | 2 | 373 |
| 54 | 5 | 5 | 58 | 3 | 3 |
| 59 | 2 | 61 | 61 | 2 | 307 |
| 62 | 3 | 3 | 67 | 2 | 337 |
| 71 | 7 | 7 | 73 | 5 | 5 |
| 74 | 5 | 5 | 79 | 3 | 3 |
| 81 | 2 | 83 | 82 | 7 | 7 |
| 83 | 2 | 251 | 86 | 3 | 3 |
| 89 | 3 | 3 | 94 | 5 | 5 |
| 97 | 5 | 5 | 98 | 3 | 3 |

Now, one notices that the case 2 is a primitive root mod $m$ occurs only if $m$ is odd. Indeed, if $(p, m) \neq 1$ then $p$ cannot be a primitive root mod $m$. In that case the value for $p$ obtained by the algorithm is not the least possible value greater than 2, but rather the least prime that is congruent to the least primitive root. But notice, if 2 is primitive for $m = q^n$, it is not primitive for $m = 2q^n$. In that case we look at the least primitive root for $m = 2q^n$. If this is prime, it is an inert prime in $K_m$ for $m = q^n$. We note however, the method of looking at $2q^n$ will not work in general if the least primitive root mod $q^n$ is composite. For example if $m = 409$ then $r = 21$, but for $m = 818$ we have $r$ is also 21. Also for $m = 271$ we have $r = 6$ and for $m = 542$ $r = 15$. But notice, for $m = 41$ we have $r = 6$ and for $m = 82$ we have $r = 7$, so there the method did work as obviously 7 is the smallest prime greater than 6.

So we need to answer, when does this method work and insure one has found the least prime primitive root that is greater than 2?

**Proposition 2.3.** *Suppose $m = q^n$ with $q > 2$ prime. Also suppose $2 \nmid r$. Then if $r$ is a least primitive root mod $2m$ it is the least primitive root mod $m$ that is greater than 2.*

*Proof.* We have $r^{\phi(2m)} \equiv 1 \bmod 2m$ and that $\phi(2m)$ is the least such exponent for which this holds and $r$ is the least number greater than 2 for which this occurs. Then, since $(2, m) = 1$ we have $\phi(2m) = \phi(2)\phi(m) = \phi(m)$. So we have $r^{\phi(m)} \equiv 1 \bmod 2m$. So, $2m|(r^{\phi(m)} - 1)$ and hence $2|(r^{\phi(m)} - 1)$ and $m|(r^{\phi(m)} - 1)$. Notice that $\phi(m)$ and $r$ are minimal with respect to this property for $r > 2$. But then $r^{\phi(m)} \equiv 1$ mod $m$ and $r$ is the lest primitive root mod $m$ which is greater than 2. $\square$

Hence the least primitive root mod $2m$ is the least primitive root mod $m$ which is greater than 2. So as long as the least primitive root mod $2m$ is prime we have, in fact, found the least inert prime. We leave open how to accomplish this if $r > 2$ and $r$ is composite. However, as among our cases this occurs only for $m = 41$, we may shorten our list and thereby obtain the following table showing the least inert prime $p$ in $K_m$ (with $m = 4$ or $q^n$) that is larger than 2. This modified table is presented next. Note that we are giving the $m$ with $3 \leq m < 100$ and $m$ not congruent to 2 mod 4 so that there is a prime $p > 2$ so that $\mathcal{O}_m/p\mathcal{O}_m$ is a finite field and in the process giving the least such $p$.

Table 2.4

| $m$ | $p$ | $m$ | $p$ |
|---|---|---|---|
| 3 | 5 | 41 | 7 |
| 4 | 3 | 43 | 3 |
| 5 | 3 | 47 | 5 |
| 7 | 3 | 49 | 3 |
| 9 | 5 | 53 | 3 |
| 11 | 7 | 59 | 11 |
| 13 | 7 | 61 | 7 |
| 17 | 3 | 67 | 7 |
| 19 | 3 | 71 | 7 |
| 23 | 5 | 73 | 5 |
| 25 | 3 | 79 | 3 |
| 27 | 5 | 81 | 5 |
| 29 | 3 | 83 | 5 |
| 31 | 3 | 89 | 3 |
| 37 | 5 | 97 | 5 |

## Section 3: The kernel group $D(\Lambda_{m,p})$

In this section we show that if $p > 2$ is unramified in $K_m/\mathbb{Q}$ then $T(\Lambda_{m,p}) =$

$D(\Lambda_{m,p})$ and is isomorphic to the cokernel of the map $h$ in the Mayer-Vietoris sequence of Reiner-Ullom. This is Theorem 3.1. Our proof also yields that this group is a quotient of the group of units of a finite field when $p$ is inert. Our result is closely related to three known results. First, [5, Corollary 1.2] yields that for any number field $K$ if $p$ is unramified in $K$ then $T(\Lambda) = D(\Lambda)$. Her proof does not use the exact sequence of Reiner-Ullom. This result was obtained as [12, Theorem 2.1] using the exact sequence of Reiner-Ullom. Last, that $D(\Lambda_{m,p})$ is isomorphic to the quotient we obtain, can be proved as a corollary to Rim's Theorem (see [13]) again using the fiber product argument. This is the content of [3, 50.45]. The advantage of our approach is that it combines these two facts and does so in a way that we obtain, in Corollary 3.2, results about $R(\Lambda_{m,p})$.

**Theorem 3.1.** *Suppose $p > 2$ is prime and that $p \nmid m$. Then $T(\Lambda_{m,p}) = D(\Lambda_{m,p}) \cong$ $cok(h)$ where $h$ is the usual map in the exact sequence of Riener-Ullom. Hence:*

$$T(\Lambda_{m,p}) = D(\Lambda_{m,p}) \cong cok(h) = (\mathbb{Z}[\zeta_{mp}]/(1-\zeta_p)\mathbb{Z}[\zeta_{mp}])^*/h(\mathbb{Z}[\zeta_{mp}]^*).$$

*Proof.* Note first that $\mathbb{Z}/p\mathbb{Z} \cong \mathbb{Z}[\zeta_p]/(1-\zeta_p)\mathbb{Z}[\zeta_p]$. Now consider the fiber product:

$$\begin{array}{ccc} \mathbb{Z}[C_p] & \longrightarrow & \mathbb{Z}[\zeta_p] \\ {\scriptstyle g}\downarrow & & \downarrow \\ \mathbb{Z} & \longrightarrow & \mathbb{Z}/p\mathbb{Z} \cong \mathbb{Z}[\zeta_p]/(1-\zeta_p)\mathbb{Z}[\zeta_p], \end{array}$$

where $f(x) = \zeta_p$ for $x$ a fixed generator of $C_p$ and $g(x) = 1$.

Now note that if $p \nmid m$ then $\mathbb{Z}[\zeta_m] \otimes \mathbb{Z}[\zeta_p] \cong \mathbb{Z}[\zeta_{mp}]$. So we apply the exact functor $\mathbb{Z}[\zeta_m] \otimes -$ to this fiber product above to obtain:

$$
\begin{array}{ccc}
\mathbb{Z}[\zeta_m]C_p & \longrightarrow & \mathbb{Z}[\zeta_{mp}] \\
\downarrow & & \downarrow \\
\mathbb{Z}[\zeta_m] & \longrightarrow & \mathbb{Z}[\zeta_m]/p\mathbb{Z}[\zeta_m] \cong \mathbb{Z}[\zeta_{mp}]/(1-\zeta_p)\mathbb{Z}[\zeta_{mp}]
\end{array}
$$

.

This fiber product yields the following exact Mayer-Vietoris sequence: (taking into account that $D(\mathbb{Z}[\zeta_{mp}])$ and $D(\mathbb{Z}[\zeta_m])$ are trivial)

$$1 \to \mathbb{Z}[\zeta_m][C_p]^* \to \mathbb{Z}[\zeta_m]^* \times \mathbb{Z}[\zeta_{mp}]^* \xrightarrow{h} (\mathbb{Z}[\zeta_{mp}]/(1-\zeta_p)\mathbb{Z}[\zeta_{mp}])^* \xrightarrow{\delta} D(\mathbb{Z}[\zeta_{mp}]C_p) \longrightarrow 0.$$

Hence the map $\delta$ is surjective and $T(\mathbb{Z}[\zeta_m]C_p) = D(\mathbb{Z}[\zeta_m]C_p) \cong cok(h)$. Now, $cok(h) = (\mathbb{Z}[\zeta_{mp}]/(1-\zeta_p)\mathbb{Z}[\zeta_{mp}])^*/h(\mathbb{Z}[\zeta_m]^* \times \mathbb{Z}[\zeta_{mp}]^*)$ and recall that $h[(u,v)] = \overline{uv}^{-1}$. However, as $\mathbb{Z}[\zeta_m]^* \subset \mathbb{Z}[\zeta_{mp}]^*$ then, by a slight abuse of notation, we can say that $h(\mathbb{Z}[\zeta_m]^* \times \mathbb{Z}[\zeta_{mp}]^*) = h(\mathbb{Z}[\zeta_{mp}]^*)$, where $h$ is now reduction modulo $(1-\zeta_p)$. □

The Theorem has two corollaries which we prove.

**Corollary 3.2, (cf [12, Theorem 4.2]).** *Let* $\Lambda_{m,p} = \mathbb{Z}[\zeta_m]C_p$ *with* $(m,p) = 1$, *then* $T^{\frac{p-1}{2}}(\Lambda_{m,p}) \subseteq R(\Lambda_{m,p}) \cap D(\Lambda_{m,p}) \subseteq T(\Lambda_{m,p})$. *Hence, if the exponent of* $T(\Lambda_{m,p})$ *is co-prime to* $\frac{p-1}{2}$, *then* $cok(h) \cong T(\Lambda_{m,p}) = R(\Lambda_{m,p}) \cap D(\Lambda_{m,p})$.

*Proof.* From Proposition 1.1 we have $T^{\frac{p-1}{2}}(\Lambda_{m,p}) \subseteq R(\Lambda_{m,p}) \cap D(\Lambda_{m,p})$. From Theorem 3.1 we have, if $(m,p) = 1$ then $cok(h) \cong T(\Lambda_{m,p}) = D(\Lambda_{m,p})$. Thus $T^{\frac{p-1}{2}}(\Lambda_{m,p}) \subseteq R(\Lambda_{m,p}) \cap D(\Lambda_{m,p}) \subseteq D(\Lambda_{m,p}) = T(\Lambda_{m,p}) \cong cok(h)$. Finally, if $T(\Lambda_{m,p})$ has exponent co-prime to $\frac{p-1}{2}$ then $T^{\frac{p-1}{2}}(\Lambda_{m,p}) = T(\Lambda_{m,p})$ and the result follows. □

**Corollary 3.3.** *If the rational prime $p > 2$ is prime in $\mathbb{Z}[\zeta_m]$ then $D(\Lambda_{m,p})$ is isomorphic to a quotient group of the group of units of the finite field $\mathbb{Z}[\zeta_m]/p\mathbb{Z}[\zeta_m]$.*

*Proof.* By Theorem 3.1 we have $T(\Lambda_{m,p}) = D(\Lambda_{m,p}) \cong cok(h) = (\mathbb{Z}[\zeta_{mp}]/(1 - \zeta_p)\mathbb{Z}[\zeta_{mp}])^*/h(\mathbb{Z}[\zeta_{mp}]^*)$. We noted in the proof that $\mathbb{Z}[\zeta_{mp}]/(1-\zeta_p)\mathbb{Z}[\zeta_{mp}] \cong \mathbb{Z}[\zeta_m]/p\mathbb{Z}[\zeta_m]$. So if $p$ is prime then $p\mathbb{Z}[\zeta_m]$ is a maximal ideal and $\mathbb{Z}[\zeta_m]/p\mathbb{Z}[\zeta_m] \cong \mathbb{F}_{p^{\phi(m)}}$, the field of $p^{\phi(m)}$ elements. □

## Section 4: Cyclotomic units and the computation of the Swan subgroup

We wish to show to what extent $T(\Lambda_{m,p})$ for $p$ prime in $\mathcal{O}_m$ is computable. By Theorem 3.1 we must compute the order of the cyclic group $cok(h) = (\mathbb{Z}[\zeta_{mp}]/(1 - \zeta_p)\mathbb{Z}[\zeta_{mp}])^*/h(\mathbb{Z}[\zeta_{mp}]^*)$ when $p$ is inert in $\mathcal{O}_m$. To begin with, let us briefly diagram the situation.

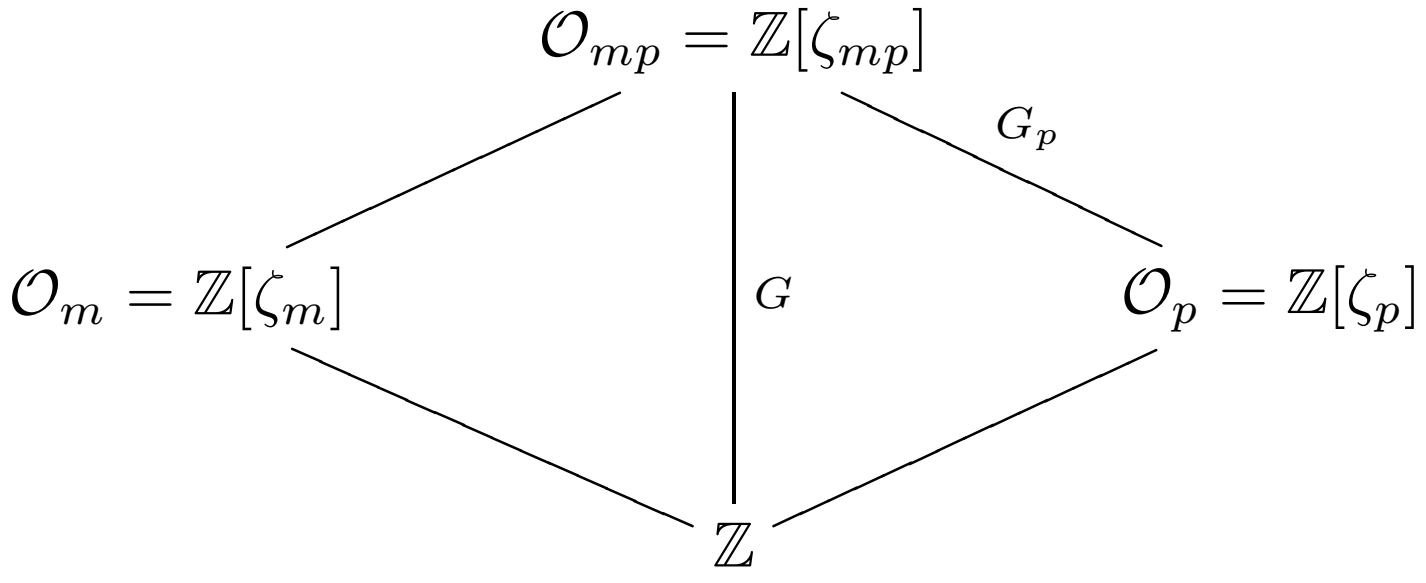

where $G = Gal(\mathbb{Q}(\zeta_{mp})/\mathbb{Q})$ and $G_p = Gal(\mathbb{Q}(\zeta_{mp})/\mathbb{Q}(\zeta_m))$. Observe that $G_p$ acts on $\mathcal{O}_{mp}$. Also, since $\mathbb{Q}(\zeta_m)$ and $\mathbb{Q}(\zeta_p)$ are linearly disjoint then, by Natural Irrationality, we have $G_p \cong Gal(\mathbb{Q}(\zeta_m)/\mathbb{Q})$.

As observed above $\mathcal{O}_{mp}/(1 - \zeta_p)\mathcal{O}_{mp} \cong \mathbb{F}_{p^{\phi(m)}}$ and, as such, the image of $h$ is contained within the unit group of this field. For computational purposes, we need the following.

**Lemma 4.1.** *Given $p$ inert in $\mathcal{O}_m$ then $\mathcal{O}_{mp}/(1-\zeta_p)\mathcal{O}_{mp} \cong \mathbb{F}_p[z]/\langle \Phi_m(z) \rangle$ where $\Phi_m(z)$ is the $m^{th}$ cyclotomic polynomial.*

*Proof.* As $p$ is inert in $\mathcal{O}_m$, the class of $\zeta_m$ in $\mathcal{O}_{mp}/(1-\zeta_p)\mathcal{O}_{mp}$ is a primitive $m^{th}$ root of unity adjoined to the (finite) field $\mathcal{O}_p/(1-\zeta_p)\mathcal{O}_p \cong \mathbb{F}_p$. $\square$

We shall denote this finite field by $\mathcal{F}_{mp}$.

When computing the image of the units in a cyclotomic field it is frequently much easier to use cyclotomic units. The immediate question is, when are the cyclotomic units the full group of units? To simplify our discussion we fix the following notation: Let the units in $\mathcal{O}_{mp}$ be denoted $E_{mp}$; denote the class number of $K_{mp}$ by $h_{mp}$; denote the class number of the maximal real subfield of $K_{mp}$ by $h^+_{mp}$; denote the cyclotomic units in $K_{mp}$ by $C_{mp}$; denote the roots of unity in $K_{mp}$ by $W_{mp}$. By [17, Theorem 4.12] we have $[E_{mp} : W_{mp}E^+_{mp}] = 2$ as $mp$ is not a prime power. Next notice that $m = 4$ or $q^n$ and that $p$ and $q$ are distinct and both greater than 2. Hence at most two primes divide $mp$. The main theorem of [15] gives that $[E^+_{mp} : C^+_{mp}] = 2^b h^+_{mp}$ where for $g$, the number of prime factors of $mp$, one defines $b$ by $b = 0$ if $g = 1$ and $b = 2^{g-2} + 1 - g$ for $g \geq 2$. Hence, in our case, (as $g = 1$ or 2) we have $[E^+_{mp} : C^+_{mp}] = h^+_{mp}$. Therefore, for our $m$ and $p$, if $h^+_{mp} = 1$ then $[E_{mp} : C_{mp}][C_{mp} : W_{mp}E^+_{mp}] = [E_{mp} : W_{mp}E^+_{mp}] = 2$ and, as $W_{mp}E^+_{mp} = W_{mp}C^+_{mp}$, and $[C_{mp} : W_{mp}C^+_{mp}] > 1$ [2, Corollary 2.3], the cyclotomic units generate the full group of units. From van der Linden [7] one finds that $h^+_{mp} = 1$ if $\phi(mp) \leq 72$. For related results on bounding class numbers see the appendix of [17] and the work of Schoof [14].

Therefore, if $\phi(mp) \leq 72$ we may just consider cyclotomic units. The question is then, for what pairs $(m,p)$ from our table does this hold?

**Remark 4.2.** *An elementary argument using properties of Euler's phi-function shows that, to require $m = q^n$ with $q > 2$ prime and $p > 2$ prime, and $\phi(mp) < 72$, forces $m \leq 37$.*

If one checks our pairs $(m,p)$ with $m \leq 37$ one finds $\phi(mp) < 72$ except for the pairs $(23,5)$ and $(37,5)$. Hence $E_{mp} = C_{mp}$ for many of the initial cases in our list. We will compute $h(\mathcal{O}^*_{mp})$ (and consequently $cok(h)$) using the generating set for $C_{mp}$ discovered by Conrad in [1]. Specifically,

**Proposition 4.3 (cf. [1]).** *Modulo roots of unity, $C_n \subseteq (\mathbb{Z}[\zeta_n])^*$ is generated by the union of the following sets of units:*

$$\{\frac{1-\zeta_d^a}{1-\zeta_d}\} \text{ } d \text{ is a prime power divisor of } n \text{ with } a \in (\mathbb{Z}/d\mathbb{Z})^*$$

$$\{1-\zeta_d^a\} \text{ } d \text{ is a \textbf{non}-prime power divisor of } n \text{ with } a \in (\mathbb{Z}/d\mathbb{Z})^*$$

*i.e.*

$$C_n = \langle \bigcup_{\substack{d|n \\ d \text{ prime power} \\ a\in(\mathbb{Z}/d\mathbb{Z})^*}} \{\frac{1-\zeta_d^a}{1-\zeta_d}\} \quad \cup \bigcup_{\substack{d|n \\ d \text{ non prime-power} \\ a\in(\mathbb{Z}/d\mathbb{Z})^*}} \{1-\zeta_d^a\}\rangle$$

*[N.B. When $n$ is a not a prime power, $C_n$ contains $\zeta_n$.]* $\square$

As such, for general $n$, one needs to consider all the possible $d$ which satisfy the above criteria. In our case $n = mp$ where $p$ is prime and $m = q, q^2$ or $q^3$ for $q$ a prime. We will consider the image of these sets of units under $h$, which amounts to considering the order of each class in $\mathcal{O}_{mp}/(1-\zeta_p) \cong \mathbb{F}_p[z]/(\Phi_m(z))$, the identification being $\zeta = \zeta_{mp} \mapsto z$ and $\zeta_p = \zeta^m$.

The goal then is to determine $|h(C_{mp})|$ using the above description. The calculations themselves were carried out using the MAPLE computer algebra system, and some procedures coded in the $C$++ language. We make a number of observations that are not only interesting from a theoretical standpoint, but were highly useful in reducing the number of calculations necessary to determine $|h(C_{mp})|$.

**Proposition 4.4.** *The homomorphism $h : \mathcal{O}^*_{mp} \to \mathcal{F}^*_{mp}$ is $G_p$-equivariant where $G_p = Gal(\mathbb{Q}(\zeta_{mp})/\mathbb{Q}(\zeta_p))$.*

*Proof.* This follows immediately from the fact that, with respect to the extension $\mathbb{Q}(\zeta_{mp})/\mathbb{Q}(\zeta_p)$, $G_p$ is the decomposition group of $1 - \zeta_p$ and the inertia group is trivial. □

In what follows, for a given $u \in \mathcal{O}^*_{mp}$, let $orb_{G_p}(u) = \{\sigma(u) : \sigma \in G_p\}$ and $orb_{G_p}(\overline{u}) = \{\sigma(h(u)) : \sigma \in G_p\}$, the conjugates of $u$ and the conjugates of $\overline{u} = h(u)$ under the respective actions of $G_p$ on $\mathcal{O}^*_{mp}$ and $\mathcal{F}^*_{mp}$. By the above observation, we see that $h(orb_{G_p}(u)) = orb_{G_p}(\overline{u})$ and that, if $v \in orb_{G_p}(u)$ then $|h(v)| = |h(u)|$.

Also, denote:

$$u_{frac,d,a} = \frac{1-\zeta_d^a}{1-\zeta_d} \text{ with } a \in (\mathbb{Z}/d\mathbb{Z})^*,\ d \text{ prime power divisor of } mp$$

$$u_{flat,d,a} = 1 - \zeta_d^a \text{ with } a \in (\mathbb{Z}/d\mathbb{Z})^*,\ d \text{ non prime power divisor of } mp$$

**Lemma 4.5.** *For each $p$ and all possible values of $m$ under consideration, $q$, $q^2$ or $q^3$, the order of $h(u_{frac,p,a})$ equals the order of $a$ mod $p$.*

*Proof.* First note that, regardless of $m$,

$$u_{frac,p,a} = \frac{1-\zeta_p^a}{1-\zeta_p} = \frac{1-\zeta^{am}}{1-\zeta^m}$$

and since, working mod $I = (1-\zeta_p)$ forces $\zeta^m \equiv 1$, then

$$\frac{1-\zeta^{am}}{1-\zeta^m} = \sum_{k=0}^{a-1} \zeta^{mk} \equiv \ a \ (mod \ I)$$

and so, $h(u_{frac,p,a}) = a \in \mathcal{F}^*_{mp}$. □

For the units denoted 'flat', it turns out that, depending on $m$, one only has to compute the order of the image of at most three of them to determine the order of all.

**Proposition 4.6.** *For each $m$, $p$ in Table 2.4, with $d$ a non prime-power divisor of $mp$,*

$$h(orb_{G_p}(1-\zeta_d)) = h(\{1-\zeta_d^a | a \in (\mathbb{Z}/d\mathbb{Z})^*\})$$

*Proof.* Since the possible values of $m$ under consideration are $q^i$ where $i \leq 3$ then $d = pq^j$ for $1 \leq j \leq i$. Now, given that $I = (1-\zeta_p) = (1-\zeta^m)$ we have $h(\zeta^a) = h(\zeta^{a'})$ if and only if $a \equiv a'$ $(mod\ m)$. This being the case, many of the units $u_{flat,d,a}$ have the same image under $h$. If $d = pq^i$ then $1-\zeta_d^a = 1-\zeta^{aq^{i-j}}$ and so $h(1-\zeta_d^a) = h(1-\zeta_d^{a'})$ if and only if $a \equiv a'$ $(mod\ q^j)$. As such, $h(\{1-\zeta_d^a | a \in (\mathbb{Z}/d\mathbb{Z})^*\})$ consists of all the 'flat' units which are distinct modulo $I$. In a similar fashion, $h(orb_{G_p}(1-\zeta^d))$ consists of the elements in the orbit of $1-\zeta$ that are unique modulo $I$. To show the equality mentioned in the statement of the proposition, we first observe the following about $G_p$.

$$\begin{aligned} G_p &= \{\sigma \in G | \sigma(\zeta_p) = \zeta_p\} \\ &\cong \{t \in (\mathbb{Z}/mp\mathbb{Z})^* | t \equiv 1 \ (mod\ p)\} \\ &\cong (\mathbb{Z}/m\mathbb{Z})^* \end{aligned}$$

So we may write, $orb_{G_p}(1-\zeta_d) = orb_{G_p}(1-\zeta^{q^{i-j}}) = \{1-\zeta^{tq^{i-j}} | t \in (\mathbb{Z}/q^i\mathbb{Z})^*\}$. As such, $h(1-\zeta^{tq^{i-j}}) = h(1-\zeta^{t'q^{i-j}})$ if and only if $t \equiv t'$ $(mod\ q^j)$. The proposition is established by observing the surjectivity of the natural maps $(\mathbb{Z}/q^i\mathbb{Z})^* \to (\mathbb{Z}/q^j\mathbb{Z})^*$ and $(\mathbb{Z}/pq^j\mathbb{Z})^* \to (\mathbb{Z}/q^j\mathbb{Z})^*$. □

**Corollary 4.7.** *For $m = q^i$ and $d = pq^j$ for $1 \le j \le i$, $|h(1-\zeta_d^a)| = |h(1-\zeta^{q^{i-j}})|$ for all $a \in (\mathbb{Z}/d\mathbb{Z})^*$.*

*Proof.* By 4.4, if $h(u) = h(\sigma(v))$ then $|h(u)| = |h(\sigma(v))| = |\sigma(h(v))| = |h(v)|$. □

As such, to compute the order (of the images) of the flat units requires the computation of only $1-\zeta_d$ for each $d = q^j$ where $m = q^i$ with $1 \le j \le i \le 3$.

The only remaining units to consider are $u_{frac,d,a}$ where $d$ is a divisor of $m$. The difficulty with these is that none are conjugate to each other under the action of $G_p$ so a technique like that in 4.6,7 does not work. Nonetheless, the list of units to check is considerably shortened by this analysis. From Table 2.4 we have for those pairs $(m,p)$ with $3 \le m < 100$ so that $m = 4$ or $q^n$ with $q > 2$ prime, $p$ is the least prime greater than 2 which is prime in $\mathbb{Z}[\zeta_m]$. From Theorem 4.4 we have, when $h^+_{mp} = 1$, and exact value for $|T(\mathbb{Z}[\zeta_m]C_p)|$. From Remark 4.2 we know this occurs for all of our pairs with $m \le 31$ except $m = 23$. Hence we have an upper bound on $T(\Lambda_{m,p})$ and know the order of $T(\Lambda_{m,p}) = D(\Lambda_{m,p})$ divides this bound. When $h^+_{mp} = 1$, we in fact know $T(\Lambda_{m,p}) = D(\Lambda_{m,p})$ is a cyclic group of the given order.

Table 4.8

| $m$ | $p$ | $\lvert T(\Lambda)\rvert$ | $m$ | $p$ | $\lvert T(\Lambda)\rvert$ |
|---|---|---|---|---|---|
| 3 | 5 | =1 | 41 | 7 | $\leq$ 973076418263561 |
| 4 | 3 | =1 | 43 | 3 | $\leq$ 121632014 |
| 5 | 3 | =1 | 47 | 5 | $\leq$ 126818393139129 |
| 7 | 3 | =2 | 49 | 3 | $\leq$ 106738298 |
| 9 | 5 | =7 | 53 | 3 | $\leq$ 23979866305 |
| 11 | 7 | =764 | 59 | 11 | $\leq$ 13443299128571962495037599194 |
| 13 | 7 | =13575 | 61 | 7 | $\leq$ 184748690907313590884125 |
| 17 | 3 | =193 | 67 | 7 | $\leq$ 57693982982891377045799212 |
| 19 | 3 | =518 | 71 | 7 | $\leq$ 2667737269476512547061391732 |
| 23 | 5 | $\leq$ 1061481 | 73 | 5 | $\leq$ 2691107610725376703967787 |
| 25 | 3 | =1181 | 79 | 3 | $\leq$ 25649083246955546 |
| 27 | 5 | =36169 | 81 | 5 | $\leq$ 45991238252616223 |
| 29 | 3 | =82465 | 83 | 5 | $\leq$ 273944187281002481299710561 |
| 31 | 3 | =231434 | 89 | 3 | $\leq$ 5532420798784332769 |
| 37 | 5 | $\leq$ 51549963049 | 97 | 5 | $\leq$ 439510970573257846930592330460696 |

We observe, in the cases where we have equality, that these orders are co-prime to $\frac{p-1}{2}$. Thus, via Corollary 3.2, we, in fact, have $cok(h) \cong T(\Lambda_{m,p}) = D(\Lambda_{m,p}) \cap R(\Lambda_{m,p})$ is a cyclic group of the given order. It is our suspicion that when one has equality (i.e. when one considers the full group of units for those $mp$ with $h^+_{mp} \neq 1$) this will still be the case.

We used the least inert prime for two reasons. First, we wished to note that such a prime exists whenever there are any inert primes. Second, in the cases $m = 3$, 4, and 5 it provides trivial Swan subgroups. This is of interest as it contrasts with the imaginary quadratic case. The Theorem [12, Theorem 4a] gives for $K = \mathbb{Q}(\sqrt{-d})$ when $d > 0$ and $d \neq 1$ or 3 that for $p$ inert in $\mathcal{O}_K$ that $C_{\frac{p+1}{2}}$ is isomorphic to a subgroup of $T(\Lambda)$. Thus the Swan subgroup in that case always nontrivial for inert primes.

However, we also note, that for $m = 2$, if we use $p = 11$ (which of course is inert) that $T(\Lambda_{m,p}) \cong C_2$. Likewise for $m = 4$ take $p = 7$ and then $T(\Lambda_{m,p}) \cong C_2$. Last, for $m = 5$ take $p = 13$ and then $T(\Lambda_{m,p}) \cong C_{17}$. Notice again, each is co-prime to $\frac{p-1}{2}$.

In light of these observations we pose the following open question: Given an algebraic number field $K$ so that there are inert primes $p$ in $K$ is it the case one may always find one so that $T(\mathcal{O}_K[C_p])$ is nontrivial and of order co-prime to $\frac{p-1}{2}$?

## SECTION 5: REFERENCES